\documentclass[10pt]{article}
\usepackage{cite}
\usepackage{mathrsfs}
\usepackage{amsfonts}
\usepackage{amsmath}
\usepackage{amsfonts,amssymb,color}
\usepackage{dsfont}
\usepackage{curves}
\usepackage{mathrsfs}
\usepackage{pifont}
\usepackage{amssymb}
\allowdisplaybreaks

\numberwithin{equation}{section}

\date{}

\def\BigRoman{\uppercase\expandafter{\romannumeral\number\count 255 }}
\def\Romannumeral{\afterassignment\BigRoman\count255=}

\begin{document}
\title{Some results on the $k$-strong parity property in a graph
}
\author{\small  Jie Wu\footnote{Corresponding
author. E-mail address: jskjdxwj@126.com}\\
\small  School of Economics and management,\\
\small  Jiangsu University of Science and Technology,\\
\small  Zhenjiang, Jiangsu 212100, China\\
}

\maketitle
\begin{abstract}
\noindent A graph $G$ has the $k$-strong parity property if for any $X\subseteq V(G)$ with $|X|$ even, $G$ contains a spanning subgraph $F$ with $d_F(u)\equiv1$ (mod 2) for each $u\in X$ and
$d_F(v)\in\{k,k+2,k+4,\ldots\}$ for each $v\in V(G)\setminus X$, where $k\geq2$ is an even integer. Kano and Matsumura proposed a characterization for a graph with the $k$-strong parity property
(M. Kano, H. Matsumura, Odd-even factors of graphs, Graphs Combin. 41 (2025) 55). In this paper, we first give a size condition for a graph to have the $k$-strong parity property. Then we establish
a signless Laplacian spectral radius condition to guarantee that a graph has the $k$-strong parity property.
\\
\begin{flushleft}
{\em Keywords:} graph; size; minimum degree; signless Laplacian spectral radius; $k$-strong parity property.

(2020) Mathematics Subject Classification: 05C50, 05C70, 90B99
\end{flushleft}
\end{abstract}

\section{Introduction}

In this paper, we only consider finite, undirected and simple graphs. Let $G$ be a graph with vertex set $V(G)$ and edge set $E(G)$. The order and size of $G$ are denoted by $|V(G)|=n$ and
$|E(G)|=e(G)$, respectively. For $v\in V(G)$, the degree of the vertex $v$ in $G$, denoted by $d_G(v)$, is the number of vertices adjacent to $v$ in $G$. Let $\delta(G)$ denote the minimum
degree of $G$. For $S\subseteq V(G)$, the subgraphs of $G$ induced by $S$ and $V(G)\setminus S$ are denoted by $G[S]$ and $G-S$, respectively. Let $S$ and $T$ be two vertex-disjoint subsets
of $V(G)$. We denote by $e_G(S,T)$ the number of edges with one end in $X$ and the other in $Y$. Given two vertex-disjoint graphs $G_1$ and $G_2$, the union of $G_1$ and $G_2$ is denoted by
$G_1\cup G_2$. The join of $G_1$ and $G_2$, denoted by $G_1\vee G_2$, is the graph obtained from $G_1\cup G_2$ by adding all edges joining a vertex of $G_1$ to a vertex of $G_2$.

Given a graph $G$ with vertex set $V(G)=\{v_1,v_2,\ldots,v_n\}$, the adjacency matrix of $G$ is defined by $A(G)=(a_{ij})_{n\times n}$, where $a_{ij}=1$ if and only if $v_iv_j\in E(G)$ and
$a_{ij}=0$ if and only if $v_iv_j\notin E(G)$. The signless Laplacian matrix of $G$ is defined by $Q(G)=D(G)+A(G)$, where $D(G)$ is the diagonal degree matrix of $G$. The largest eigenvalues
of $A(G)$ and $Q(G)$, denoted by $\rho(G)$ and $q(G)$, are called the adjacency spectral radius and the signless Laplacian spectral radius of $G$, respectively.

Let $g$ and $f$ be two nonnegative integer-valued functions defined on $V(G)$ satisfying $g(v)\leq f(v)$ for any $v\in V(G)$. Then a $(g,f)$-factor of $G$ is defined as a spanning subgraph
$F_1$ of $G$ such that $g(v)\leq d_{F_1}(v)\leq f(v)$ for all $v\in V(G)$. A $(g,f)$-parity factor of $G$ is a spanning subgraph $F_2$ with $d_{F_2}(v)\equiv g(v)\equiv f(v)$ (mod 2) and
$g(v)\leq d_{F_2}(v)\leq f(v)$ for all $v\in V(G)$. If $g(v)=1$ and $f(v)=b$ for all $v\in V(G)$, then a $(g,f)$-parity factor is simply called a $(1,b)$-odd factor, where $b$ is a positive
odd integer. Obviously, a $(1,1)$-odd factor is simply called a perfect matching.

Niessen and Randerath \cite{NR}, Enomoto and Hagita \cite{EH} proposed some sufficient conditions for graphs to possess $k$-factors. Lots of researchers \cite{LP,Dai,Wp,ZSL1,ZZS} proved
some results on the existence of a $(1,2)$-factor in a graph. Miao and Li \cite{ML}, Zhou \cite{Zs}, Zhou and Wu \cite{ZW} provided some sufficient conditions for the existence of a
$(1,b)$-factor in a graph. Many efforts have been devoted to finding sufficient conditions for graphs to possess $(a,b)$-factors by using various graphic parameters such as independence
number \cite{KL}, degree condition \cite{Zr} and Fan-type condition \cite{Mf}.

O \cite{O}, Zhou, Sun and Zhang \cite{ZSZ}, Zhou and Zhang \cite{ZZ}, Zhou, Zhang and Liu \cite{ZZL}, Zhou \cite{Zt} presented some spectral conditions for a connected graph to have a perfect
matching. Amahashi \cite{A} obtained a characterization for a graph with a $(1,b)$-odd factor. Kim, O, Park and Ree \cite{KOPR} established a connection between eigenvalue and $(1,b)$-odd
factors in graphs. Hua and Zhang \cite{HZ} provided some sufficient conditions according to the distance spectral radius which guarantee the existence of $(1,b)$-odd factors in 1-binding graphs.
Liu and Lu \cite{LL} posed a degree condition for a graph to possess an $(a,b)$-parity factor. Yang, Zhang, Lu and Liu \cite{YZLL} presented an independence number and connectivity condition for
the existence of an $(a,b)$-parity factor in a graph. Some other results on graph factors see \cite{Wa,Wc,WZ,ZSL}.

A graph $G$ has the strong parity property if for every subset $X\subseteq V(G)$ with $|X|$ even, $G$ contains a spanning subgraph $F$ with $\delta(F)\geq1$, $d_F(u)\equiv1$ (mod 2) for each
$u\in X$, and $d_F(v)\equiv0$ (mod 2) for each $v\in V(G)\setminus X$. Bujt\'as, Jendrol' and Tuza \cite{BJT} introduced the concept of strong parity property and obtained some sufficient
conditions for graphs to admit the strong parity property. Lu, Yang and Zhang \cite{LYZ} showed a necessary and sufficient condition for the existence of the strong parity property in a graph.
Hasanvand \cite{H} studied the strong parity property of edge-connected graphs. Zhou, Zhang and Bian \cite{ZZB} established a relationship between adjacency spectral radius and the strong parity
property in a connected graph. Zhang and Fan \cite{ZF} provided edge number and spectral radius conditions for the strong parity property in a graph, respectively. Wang \cite{Wr} gave two
sufficient conditions for graphs to possess the strong parity property.

Kano and Matsumura \cite{KM} generalized the strong parity property of graphs to the $k$-strong parity property. Let $k\geq2$ be an even integer. A graph $G$ has the $k$-strong parity property
if for any $X\subseteq V(G)$ with $|X|$ even, $G$ contains a spanning subgraph $F$ with $d_F(u)\equiv1$ (mod 2) for each $u\in X$ and $d_F(v)\in\{k,k+2,k+4,\ldots\}$ for each $v\in V(G)\setminus X$.
Kano and Matsumura \cite{KM} provided a characterization for a graph with the $k$-strong parity property, which is a generalization of Lu, Yang and Zhang's result.

\medskip

\noindent{\textbf{Theorem 1.1}} (Kano and Matsumura \cite{KM}). Let $k\geq2$ be an even integer, and let $G$ be a connected graph. Then $G$ has the $k$-strong parity property if and only if
$$
c(G-S)\leq\sum\limits_{v\in S}{d_{G}(v)}-k|S|+1 
$$
for every $S\subseteq V(G)$, where $c(G-S)$ is the number of components in $G-S$.

\medskip

Motivated by \cite{KM} directly, it is natural and interesting to find a sufficient condition to ensure that a graph has the $k$-strong parity property. In this paper, we first establish a
relationship between the size and the $k$-strong parity property in a connected graph. Then we put forward a signless Laplacian spectral radius condition for the $k$-strong parity property in
a connected graph.

\medskip

\noindent{\textbf{Theorem 1.2.}} Let $k\geq2$ be an even integer, and let $G$ be a connected graph with minimum degree $\delta\geq k+1$ and order
$n\geq\max\Big\{\frac{(\delta^{2}+(2-k)\delta-2k-1)^{2}}{6(\delta-k)}+\delta+k+1,\delta^{2}+4\delta-k^{2}+6\Big\}$. If
$$
e(G)\geq e(K_{\delta}\vee(K_{n-(\delta-k+1)\delta-1}\cup((\delta-k)\delta+1)K_1)),
$$
then $G$ has the $k$-strong parity property unless $G=K_{\delta}\vee(K_{n-(\delta-k+1)\delta-1}\cup((\delta-k)\delta+1)K_1)$.

\medskip

\noindent{\textbf{Theorem 1.3.}} Let $k\geq2$ be an even integer, and let $G$ be a connected graph of order $n\geq\frac{2k+3+\sqrt{8k^{2}-8k+1}}{2}$. If
$$
q(G)>2(n-2)+\frac{2k-2}{n-1},
$$
then $G$ has the $k$-strong parity property.

\medskip

\section{The proof of Theorem 1.2}

Zheng, Li, Luo and Wang \cite{ZLLW} proved the following result, which will be used in the proof of Theorem 1.2.

\medskip

\noindent{\textbf{Lemma 2.1}} (Zheng, Li, Luo and Wang \cite{ZLLW}). Let $\sum\limits_{i=1}^{t}n_i=n-s$. If $n_1\geq n_2\geq\cdots\geq n_t\geq p\geq1$ and $n_1<n-s-p(t-1)$, then
$$
e(K_s\vee(K_{n_1}\cup K_{n_2}\cup\cdots\cup K_{n_t}))<e(K_s\vee(K_{n-s-p(t-1)}\cup(t-1)K_p)).
$$

\medskip

In what follows, we put forward a proof of Theorem 1.2, which shows a size condition for the $k$-strong parity property in a connected graph.

\medskip

\noindent{\it Proof of Theorem 1.2.} Suppose, to the contrary, that a connected graph $G$ does not satisfy the $k$-strong parity property. Then it follows from Theorem 1.1 that
\begin{align}\label{eq:2.1}
c(G-S)\geq\sum\limits_{v\in S}{d_G(v)}-k|S|+2\geq(\delta-k)|S|+2
\end{align}
for some nonempty subset $S$ of $V(G)$. Let $|S|=s$. Then $G$ is a spanning subgraph of $G_1=K_s\vee(K_{n_1}\cup K_{n_2}\cup\cdots\cup K_{n_{(\delta-k)s+2}})$, where $n_1,n_2,\ldots,n_{(\delta-k)s+2}$
are positive integers with $n_1\geq n_2\geq \cdots\geq n_{(\delta-k)s+2}$ and $\sum\limits_{i=1}^{(\delta-k)s+2}{n_i}=n-s$. Thus, we conclude
\begin{align}\label{eq:2.2}
e(G)\leq e(G_1),
\end{align}
where the equality holds if and only if $G=G_1$. We shall proceed by discussing the following two possible cases.

\noindent{\bf Case 1.} $s\leq\delta$.

Let $G_2=K_s\vee(K_{n-s-(\delta+1-s)((\delta-k)s+1)}\cup((\delta-k)s+1)K_{\delta+1-s})$. Note that $\delta(G)=\delta$. Thus, we possess $\delta(G_1)\geq\delta(G)=\delta$, and so $n_{(\delta-k)s+2}\geq\delta+1-s$.
According to Lemma 2.1, we get
\begin{align}\label{eq:2.3}
e(G_1)\leq e(G_2)
\end{align}
with equality if and only if $(n_1,n_2,\ldots,n_{(\delta-k)s+2})=(n-s-(\delta+1-s)((\delta-k)s+1),\delta+1-s,\ldots,\delta+1-s)$.

Let $G_*=K_{\delta}\vee(K_{n-(\delta-k+1)\delta-1}\cup((\delta-k)\delta+1)K_1)$. By a direct computation, we obtain
$$
e(G_*)=\binom{n-(\delta-k)\delta-1}{2}+\delta((\delta-k)\delta+1)
$$
and
$$
e(G_2)=\binom{n-(\delta+1-s)((\delta-k)s+1)}{2}+s((\delta-k)s+1)(\delta+1-s)+((\delta-k)s+1)\binom{\delta+1-s}{2}.
$$
Furthermore, we get
\begin{align}\label{eq:2.4}
e(G_*)-e(G_2)=\frac{1}{2}(\delta-s)f(s),
\end{align}
where $f(s)=(\delta-k)^{2}s^{3}-(\delta-k)(\delta^{2}+(2-k)\delta-2k-1)s^{2}+(\delta-k)(2n-2\delta-k-4)s-2(\delta-k-1)n+\delta^{3}+(2-2k)\delta^{2}+(k^{2}-2k+1)\delta-3k-2$, $s\in[1,\delta]$. It is straightforward
to check that the derivative function of $f(x)$ is
$$
f'(x)=3(\delta-k)^{2}x^{2}-2(\delta-k)(\delta^{2}+(2-k)\delta-2k-1)x+(\delta-k)(2n-2\delta-k-4).
$$
Notice that the symmetry axis of the parabola $f'(x)$ is $x=\frac{\delta^{2}+(2-k)\delta-2k-1}{3(\delta-k)}$. Then we possess
\begin{align*}
f'(x)\geq&f'\left(\frac{\delta^{2}+(2-k)\delta-2k-1}{3(\delta-k)}\right)\\
=&3(\delta-k)^{2}\left(\frac{\delta^{2}+(2-k)\delta-2k-1}{3(\delta-k)}\right)^{2}\\
&-2(\delta-k)(\delta^{2}+(2-k)\delta-2k-1)\left(\frac{\delta^{2}+(2-k)\delta-2k-1}{3(\delta-k)}\right)\\
&+(\delta-k)(2n-2\delta-k-4)\\
=&(\delta-k)(2n-2\delta-k-4)-\frac{(\delta^{2}+(2-k)\delta-2k-1)^{2}}{3}\\
\geq&(\delta-k)\left(\frac{(\delta^{2}+(2-k)\delta-2k-1)^{2}}{3(\delta-k)}+2\delta+2k+2-2\delta-k-4\right)\\
&-\frac{(\delta^{2}+(2-k)\delta-2k-1)^{2}}{3}\\
=&(\delta-k)(k-2)\\
\geq&0
\end{align*}
by $k\geq2$, $\delta\geq k+1$ and $n\geq\frac{(\delta^{2}+(2-k)\delta-2k-1)^{2}}{6(\delta-k)}+\delta+k+1$, which yields that $f(x)$ is increasing in the interval $[1,\delta]$. Hence, we obtain
\begin{align}\label{eq:2.5}
f(s)\geq&f(1)\nonumber\\
=&2n-\delta^{2}+(k-2)\delta-2\nonumber\\
\geq&\frac{(\delta^{2}+(2-k)\delta-2k-1)^{2}}{3(\delta-k)}+2\delta+2k+2-\delta^{2}+(k-2)\delta-2\nonumber\\
=&\frac{(\delta^{2}+(2-k)\delta-2k-1)^{2}}{3(\delta-k)}-\delta^{2}+k\delta+2k\nonumber\\
\geq&\frac{((k+1)\delta+(2-k)\delta-2k-1)(\delta^{2}+(2-k)\delta-2k-1)}{3(\delta-k)}-\delta^{2}+k\delta+2k\nonumber\\
=&\frac{(3\delta-2k-1)(\delta^{2}+(2-k)\delta-2k-1)}{3(\delta-k)}-\delta^{2}+k\delta+2k\nonumber\\
>&\frac{(3\delta-3k)(\delta^{2}+(2-k)\delta-2k-1)}{3(\delta-k)}-\delta^{2}+k\delta+2k\nonumber\\
=&2\delta-1\nonumber\\
>&0
\end{align}
by $s\in[1,\delta]$, $k\geq2$, $\delta\geq k+1$ and $n\geq\frac{(\delta^{2}+(2-k)\delta-2k-1)^{2}}{6(\delta-k)}+\delta+k+1$. In terms of \eqref{eq:2.4}, \eqref{eq:2.5} and $s\leq\delta$, we conclude
\begin{align}\label{eq:2.6}
e(G_2)\leq e(G_*)=e(K_{\delta}\vee(K_{n-(\delta-k+1)\delta-1}\cup((\delta-k)\delta+1)K_1)),
\end{align}
where the first equality holds if and only if $G_2=G_*$.

It follows from \eqref{eq:2.2}, \eqref{eq:2.3} and \eqref{eq:2.6} that
$$
e(G)\leq e(K_{\delta}\vee(K_{n-(\delta-k+1)\delta-1}\cup((\delta-k)\delta+1)K_1)),
$$
with equality if and only if $G=K_{\delta}\vee(K_{n-(\delta-k+1)\delta-1}\cup((\delta-k)\delta+1)K_1)$, a contradiction.

\noindent{\bf Case 2.} $s\geq\delta+1$.

Recall that $G_1=K_s\vee(K_{n_1}\cup K_{n_2}\cup\cdots\cup K_{n_{(\delta-k)s+2}})$. Let $G_3=K_s\vee(K_{n-(\delta-k+1)s-1}\cup((\delta-k)s+1)K_1)$. Using Lemma 2.1, we obtain
\begin{align}\label{eq:2.7}
e(G_1)\leq e(G_3),
\end{align}
with equality holding if and only if $(n_1,n_2,\ldots,n_{(\delta-k)s+2})=(n-(\delta-k+1)s-1,1,\ldots,1)$. Recall that $G_*=K_{\delta}\vee(K_{n-(\delta-k+1)\delta-1}\cup((\delta-k)\delta+1)K_1)$. By
a direct calculation, we possess
\begin{align}\label{eq:2.8}
e(G_*)-e(G_3)=&\binom{n-(\delta-k)\delta-1}{2}+\delta((\delta-k)\delta+1)-\binom{n-(\delta-k)s-1}{2}-s((\delta-k)s+1)\nonumber\\
=&\frac{1}{2}(s-\delta)g(s),
\end{align}
where $g(s)=2(\delta-k)n-(\delta-k)(\delta-k+2)s-\delta^{3}+(2k-2)\delta^{2}-(k^{2}-2k+3)\delta+3k-2$.

\noindent{\bf Subcase 2.1.} $s\geq\delta+2k+4$.

Notice that $n\geq(\delta-k+1)s+2$. Combining this with $\delta\geq k+1$ and $s\geq\delta+2k+4$, we get
\begin{align}\label{eq:2.9}
g(s)\geq&2(\delta-k)((\delta-k+1)s+2)-(\delta-k)(\delta-k+2)s\nonumber\\
&-\delta^{3}+(2k-2)\delta^{2}-(k^{2}-2k+3)\delta+3k-2\nonumber\\
=&(\delta-k)^{2}s-\delta^{3}+(2k-2)\delta^{2}-(k^{2}-2k-1)\delta-k-2\nonumber\\
\geq&(\delta-k)^{2}(\delta+2k+4)-\delta^{3}+(2k-2)\delta^{2}-(k^{2}-2k-1)\delta-k-2\nonumber\\
=&(2k+2)\delta^{2}-(4k^{2}+6k-1)\delta+2k^{3}+4k^{2}-k-2.
\end{align}
Note that
$$
\frac{4k^{2}+6k-1}{2(2k+2)}<k+1\leq\delta.
$$
Together with \eqref{eq:2.9}, we deduce
\begin{align}\label{eq:2.10}
g(s)\geq&(2k+2)\delta^{2}-(4k^{2}+6k-1)\delta+2k^{3}+4k^{2}-k-2\nonumber\\
=&(2k+2)(k+1)^{2}-(4k^{2}+6k-1)(k+1)+2k^{3}+4k^{2}-k-2\nonumber\\
=&1.
\end{align}

According to \eqref{eq:2.8}, \eqref{eq:2.10} and $s\geq\delta+2k+4$, we conclude
$$
e(G_3)<e(G_*)=e(K_{\delta}\vee(K_{n-(\delta-k+1)\delta-1}\cup((\delta-k)\delta+1)K_1)).
$$
Combining this with \eqref{eq:2.2} and \eqref{eq:2.7}, we obtain
$$
e(G)\leq e(G_1)\leq e(G_3)<e(K_{\delta}\vee(K_{n-(\delta-k+1)\delta-1}\cup((\delta-k)\delta+1)K_1)),
$$
which contradicts $e(G)\geq e(K_{\delta}\vee(K_{n-(\delta-k+1)\delta-1}\cup((\delta-k)\delta+1)K_1))$.

\noindent{\bf Subcase 2.2.} $\delta+1\leq s\leq\delta+2k+3$.

By virtue of $\delta\geq k+1$, $n\geq\delta^{2}+4\delta-k^{2}+6$ and $\delta+1\leq s\leq\delta+2k+3$, we have
\begin{align*}
g(s)=&2(\delta-k)n-(\delta-k)(\delta-k+2)s-\delta^{3}+(2k-2)\delta^{2}-(k^{2}-2k+3)\delta+3k-2\\
\geq&2(\delta-k)(\delta^{2}+4\delta-k^{2}+6)-(\delta-k)(\delta-k+2)(\delta+2k+3)\\
&-\delta^{3}+(2k-2)\delta^{2}-(k^{2}-2k+3)\delta+3k-2\\
=&\delta^{2}-(2k-3)\delta+k^{2}-3k-2\\
\geq&(k+1)^{2}-(2k-3)(k+1)+k^{2}-3k-2\\
=&2.
\end{align*}
Combining this with \eqref{eq:2.8} and $s\geq\delta+1$, we possess
\begin{align}\label{eq:2.11}
e(G_3)<e(G_*)=e(K_{\delta}\vee(K_{n-(\delta-k+1)\delta-1}\cup((\delta-k)\delta+1)K_1)).
\end{align}

It follows from \eqref{eq:2.2}, \eqref{eq:2.7} and \eqref{eq:2.11} that
$$
e(G)\leq e(G_1)\leq e(G_3)<e(K_{\delta}\vee(K_{n-(\delta-k+1)\delta-1}\cup((\delta-k)\delta+1)K_1)),
$$
which contradicts $e(G)\geq e(K_{\delta}\vee(K_{n-(\delta-k+1)\delta-1}\cup((\delta-k)\delta+1)K_1))$. This completes the proof of Theorem 1.2. \hfill $\Box$

\section{The proof of Theorem 1.3}

We first introduce a lemma, which will be used in the proof of Theorem 1.3.

\medskip

\noindent{\textbf{Lemma 3.1}} (Das \cite{Dm}). Let $G$ be a graph of order $n$. Then
$$
q(G)\leq\frac{2e(G)}{n-1}+n-2.
$$

\medskip

In what follows, we verify Theorem 1.3, which establishes a relationship between signless Laplacian spectral radius and the $k$-strong parity property in a connected graph.

\medskip

\noindent{\it Proof of Theorem 1.3.} Suppose, to the contrary, that a connected graph $G$ does not satisfy the $k$-strong parity property. In terms of Theorem 1.1, we conclude
\begin{align}\label{eq:3.1}
c(G-S)\geq\sum\limits_{v\in S}{d_G(v)}-k|S|+2
\end{align}
for some nonempty subset $S$ of $V(G)$. Let $|S|=s$ and $c(G-S)=t$, and let $C_1,C_2,\ldots,C_t$ be the components of $G-S$. Then $n=|S|+|V(C_1)|+\cdots+|V(C_t)|\geq s+t$.

\noindent{\bf Claim 1.} $t\geq1$.

\noindent{\it Proof.} Assume that $t=0$. Then $s=n$. By means of \eqref{eq:3.1}, we obtain
$$
e(G)=e(G[S])=\frac{1}{2}\sum\limits_{v\in S}{d_G(v)}\leq\frac{1}{2}(kn-2).
$$
Combining this with Lemma 3.1 and $n\geq\frac{2k+3+\sqrt{8k^{2}-8k+1}}{2}$, we deduce
$$
q(G)\leq\frac{2e(G)}{n-1}+n-2\leq\frac{kn-2}{n-1}+n-2<2(n-2)+\frac{2k-2}{n-1},
$$
which is a contradiction to $q(G)>2(n-2)+\frac{2k-2}{n-1}$. Claim 1 is proved. \hfill $\Box$

Notice that $|V(C_i)|\geq1$ by Claim 1, where $1\leq i\leq t$. Then it follows from \eqref{eq:3.1} that
\begin{align}\label{eq:3.2}
e(G)=&e(G[S])+\sum_{i=1}^{t}e_G(S,V(C_i))+\sum_{i=1}^{t}e(C_i)\nonumber\\
=&\sum\limits_{v\in S}{d_G(v)}-\frac{1}{2}\sum\limits_{v\in S}{d_{G[S]}(v)}+\sum_{i=1}^{t}e(C_i)\nonumber\\
\leq&\sum\limits_{v\in S}{d_G(v)}+\sum_{i=1}^{t}e(C_i)\nonumber\\
\leq&\sum\limits_{v\in S}{d_G(v)}+e(K_{n-s-t+1}\cup(t-1)K_1)\nonumber\\
\leq&t+ks-2+\binom{n-s-t+1}{2}\nonumber\\
=&\frac{1}{2}M,
\end{align}
where $M=(n-s-t)(n-s-t+1)+2t+2ks-4$. Recall that $n\geq s+t$. Then we shall proceed by discussing the following two possible cases.

\noindent{\bf Case 1.} $s+t\leq n\leq s+t+k-1$.

Let $\psi(n)=M=n^{2}-(2s+2t-1)n+s^{2}+(2t+2k-1)s+t^{2}+t-4$. Then the symmetry axis of $\psi(n)$ is $n=s+t-\frac{1}{2}$, which implies that $\psi(n)$ is increasing in the interval $[s+t,s+t+k-1]$.
Combining this with Claim 1 and $s+t\leq n\leq s+t+k-1$, we conclude
\begin{align}\label{eq:3.3}
M=&\psi(n)\leq\psi(s+t+k-1)\nonumber\\
=&k(k-1)+2k(s+t)-(2k-2)t-4\nonumber\\
\leq&k(k-1)+2kn-2k-2\nonumber\\
=&2kn+k(k-3)-2.
\end{align}

According to \eqref{eq:3.2}, \eqref{eq:3.3}, Lemma 3.1 and $n\geq\frac{2k+3+\sqrt{8k^{2}-8k+1}}{2}$, we have
\begin{align*}
q(G)\leq&\frac{2e(G)}{n-1}+n-2\\
\leq&\frac{M}{n-1}+n-2\\
\leq&\frac{2kn+k(k-3)-2}{n-1}+n-2\\
\leq&2(n-2)+\frac{2k-2}{n-1},
\end{align*}
which contradicts $q(G)>2(n-2)+\frac{2k-2}{n-1}$.

\noindent{\bf Case 2.} $n\geq s+t+k$.

Let $\varphi(s)=M=s^{2}-(2n-2t-2k+1)s+(n-t)(n-t+1)+2t-4$. It is straightforward to check that the derivative function of $\varphi(s)$ is
$$
\varphi'(s)=2s-2n+2t+2k-1.
$$
According to $n\geq s+t+k$, we obtain
$$
\varphi'(s)=2s-2n+2t+2k-1\leq2s-2(s+t+k)+2t+2k-1=-1<0,
$$
which yields that $\varphi(s)$ is decreasing in the interval $[1,+\infty)$. Thus, we possess
\begin{align}\label{eq:3.4}
M=&\varphi(s)\leq\varphi(1)\nonumber\\
=&n^{2}-3n+2k-(t-1)(2n-t-4)\nonumber\\
\leq&n^{2}-3n+2k,
\end{align}
where the first inequality holds from the fact $s\geq1$ and the last inequality holds from the fact $t\geq1$ and $n\geq s+t+k\geq t+3$. In terms of \eqref{eq:3.2}, \eqref{eq:3.4} and Lemma 3.1, we get
\begin{align*}
q(G)\leq&\frac{2e(G)}{n-1}+n-2\\
\leq&\frac{M}{n-1}+n-2\\
\leq&\frac{n^{2}-3n+2k}{n-1}+n-2\\
=&2(n-2)+\frac{2k-2}{n-1},
\end{align*}
which contradicts $q(G)>2(n-2)+\frac{2k-2}{n-1}$. This completes the proof of Theorem 1.3. \hfill $\Box$

\section{Concluding remarks}

Kano and Matsumura \cite{KM} generalized the strong parity property of graphs to the $k$-strong parity property. Furthermore, they also provided a necessary and sufficient condition for a graph to have the
$k$-strong parity property, which is a generalization of Lu, Yang and Zhang's result \cite{LYZ}. In this paper, we also investigate the $k$-strong parity property of a graph and obtain two sufficient conditions
for a graph to have the $k$-strong parity property via the size and the signless Laplacian spectral radius, respectively. Indeed, there are very few results on the $k$-strong parity property of a graph. Hence,
it is natural and interesting to put forward some new sufficient conditions to ensure that a graph has the $k$-strong parity property.

\section*{Declaration of competing interest}

\medskip

The author declares that he has no known competing financial interests or personal relationships that could have appeared to influence the work
reported in this paper.

\section*{Data availability}

\medskip

No data was used for the research described in the article.

\medskip


\end{document}